\documentclass[12pt, a4paper, reqno]{amsart}
\usepackage{inputenc}

\usepackage[all]{xy}

\usepackage{booktabs}
\usepackage{tabularx}
\usepackage{geometry}

\usepackage{fancyhdr}
\fancyheadoffset{\dimexpr\oddsidemargin-1in} 

\usepackage{setspace}

\usepackage{tikz}
\usepackage[normalem]{ulem}
\usepackage{soul}
\usepackage{color}
\usepackage{upgreek}

\setstcolor{red}
\makeatletter
\@namedef{subjclassname@2010}{%
	\textup{2010} Mathematics Subject Classification}
\makeatother

\makeatletter
\newcommand{\biggg}{\bBigg@{3}}

\newcommand{\Biggg}{\bBigg@{3.5}}

\newcommand{\bigggg}{\bBigg@{4}}

\newcommand{\Bigggg}{\bBigg@{4.5}}

\makeatother

\numberwithin{equation}{section}

\usepackage{verbatim}

\allowdisplaybreaks[3]

\usepackage{amsmath}
\usepackage{amsthm}
\usepackage{amssymb}
\usepackage{mathtools}
\usepackage{mathrsfs}
\usepackage{amsxtra}

\newtheorem{theorem1}{Theorem}[section]
\newtheorem{theorem2}[theorem1]{Theorem}
\newtheorem{lemma1}{Lemma}[section]
\newtheorem{lemma2}[lemma1]{Lemma}
\newtheorem{lemma3}[lemma1]{Lemma}
\newtheorem{lemma4}[lemma1]{Lemma}
\newtheorem{lemma5}[lemma1]{Lemma}
\newtheorem{lemma6}[lemma1]{Lemma}
\newtheorem{lemma7}[lemma1]{Lemma}
\newtheorem{lemma8}[lemma1]{Lemma}
\newtheorem{lemma9}[lemma1]{Lemma}
\newtheorem{lemma10}[lemma1]{Lemma}
\newtheorem{lemma11}[lemma1]{Lemma}
\newtheorem{lemma12}[lemma1]{Lemma}
\newtheorem{lemma13}[lemma1]{Lemma}

\begin{document}
\vglue -5mm

\title[On the error term concerning the number of subgroups of 
 $\mathbb{Z}_l\times\mathbb{Z}_m\times\mathbb{Z}_n$ with $lmn\leqslant x$]
{On the error term concerning the number of subgroups of  
	$\mathbb{Z}_l\times\mathbb{Z}_m\times\mathbb{Z}_n$ with $lmn\leqslant x$}


\author{Li JiaMin, Jing Ma, Jia Zhang}

\address{%
	JiaMin Li
	\\
	School of Mathematics
	\\
	Jilin University
	\\
	Changchun 130012
	\\
	P. R. China
}
\email{lijiamin@mail.sdu.edu.cn}

\address{%
	Jing Ma
	\\
	School of Mathematics
	\\
	Jilin University
	\\
	Changchun 130012
	\\
	P. R. China
}
\email{jma@jlu.edu.cn}

\address{
	Jia Zhang
	\\
	School of Mathematics
	\\
	Jilin University
	\\
	Changchun 130012
	\\
	P. R. China
}
\email{zhagnjia22@mails.jlu.edu.cn}

\date{\today}

\begin{abstract}
Let $Z_n$ denote the additive group of residue classes modulo $n$. 
Let $c(l,m,n)$ denote the 
number of cyclic subgroups of $Z_l\times Z_m\times Z_n$. For any $x\geqslant 1$, we 
consider the asymptotic behavior of $D_{3c}(x):=\sum_{lmn\leqslant x}c(l,m,n)$, obtain an asymptotic formula by complex method, 
and get an upper bound for the integral mean-square of the error term in that asymptotic formula. 
\end{abstract}

\maketitle

\section{Inteoduction}

Let $\mathbb{Z}_n$ be the additive group of residue classes modulo $n$. For arbitrary positive integers $m$ 
and 
$n$, the group $G:=\mathbb{Z}_m\times\mathbb{Z}_n$ is isomorphic to 
$\mathbb{Z}_{\rm{gcd}(m,n)}\times
{\mathbb{Z}_{\rm{lcm}(m,n)}}$. When $\rm{gcd}(m,n)=1$, $G$ is cyclic and isomorphic to $\mathbb{Z}_{mn}$. 
When $\rm{gcd}(m,n)>1$, $G$ is said to have rank two.

Let $s(n_1, \dots, n_r)$ and $c(n_1, \dots, n_r)$ denote the number of subgroups 
and the number of cyclic subgroups of $\mathbb{Z}_{n_1}\times \cdots \times \mathbb{Z}_{n_r}$, respectively.
 Then 
\begin{equation}\label{s_mn}
s(m,n)=\prod_{p}s(p^{\nu_p(m)},p^{\nu_p(n)}),
\end{equation}
\begin{equation}\label{c_mn}
c(m,n)=\prod_{p}c(p^{\nu_p(m)},p^{\nu_p(n)}),
\end{equation}
for all $m,n\in\mathbb{N}$. For the rank two $p$-group $\mathbb{Z}_{p^{a}}\times
\mathbb{Z}_{p^{b}}$ with $1\leqslant a \leqslant b$, 
we have
\begin{equation}\label{spp}
	s(p^a,p^b)=\dfrac{(b-a+1)p^{a+2}-(b-a-1)p^{a+1}-(a+b+3)p+(a+b+1)}{(p-1)^2}
\end{equation}
and
\begin{equation}\label{cpp}
	c(p^a,p^b)=2(1+p+p^2+\cdots+p^{a-1})+(b-a+1)p^a.
\end{equation}
Formula \eqref{spp} was deduced by using Goursat's Lemma for groups in \cite{Calugareanu.G,Petrillo}, and using the concept of the fundamental group lattice in \cite{Tarnauceanu2007,Tarnauceanu2010}. 
Formula \eqref{cpp} was proved in \cite{Tarnauceanu2010}.
One can compute $s(m,n)$ and $c(m,n)$ by using \eqref{s_mn}, \eqref{spp} and \eqref{c_mn}, \eqref{cpp}, 
respectively. 
Hampejs ect. \cite{M_Hampejs} give  more compact identities for $m,n\in \mathbb{N}$   
\begin{equation}\label{s_mn1}
	s(m,n)=\sum_{d|m,e|n}{\rm{gcd}}(d,e)=\sum_{d|{\rm{gcd}}(m,n)}\phi(d)\tau(m/d)\tau(n/d),
\end{equation}
\begin{equation}\label{c_mn1}
	c(m,n)=\sum_{d|m,e|n}\phi({\rm{gcd}}(d,e))=\sum_{d|{\rm{gcd}}(m,n)}(\mu*\phi)(d)\tau(m/d)\tau(n/d),
\end{equation}
where $\phi$ is the Euler function.
 

Let $x\geqslant 1$ be a real number and $r\in \mathbb{N}$. Define
\begin{align*}
	&S_r^{(1)}(x):=\sum_{n_1\leqslant x,\dots,n_r \leqslant x}s(n_1,\dots,n_r),\qquad
	C_r^{(1)}(x):=\sum_{n_1\leqslant x,\dots,n_r \leqslant x}c(n_1,\dots,n_r),
\\
	&S_2^{(2)}(x):=\underset{{\rm{gcd}}(m,n)>1}{\sum_{m,n\leqslant x}}s(m,n),
	\hskip21mm
	C_2^{(2)}(x):=\underset{{\rm{gcd}}(m,n)>1}
	              {\sum_{m,n \leqslant x}}c(m,n).
\end{align*} 
Then, $S^{(1)}_r(x)$ and $C^{(1)}_r(x)$ are the total number of subgroups and the  total  number of cyclic subgroups of all $\mathbb{Z}_{n_1} \times\cdots\times\mathbb{Z}_{n_r}$ with $n_i\leqslant x$, respectively;  
$S^{(2)}_2(x)$ and $C^{(2)}_2(x)$ are the  total  number of subgroups and the total  number of cyclic subgroups 
of all rank two groups $\mathbb{Z}_m\times\mathbb{Z}_n$ with $m,n\leqslant x$ respectively.

For $r=2$, Nowak and T\'oth \cite{W_G_Nowak} proved that, for $j=1,2$,
\begin{equation}\label{WGN1}
	S^{(j)}_2(x)=x^2\sum_{l=0}^{3}A_{j,l}\log^lx+O(x^{1117/701+\varepsilon}),
\end{equation}
\begin{equation}\label{WGN2}
	C^{(j)}_2(x)=x^2\sum_{l=0}^{3}B_{j,l}\log^lx+O(x^{1117/701+\varepsilon}),
\end{equation}
where $A_{j,l}, B_{j,l}$ $(1\leqslant j\leqslant2,\ 0\leqslant l \leqslant3)$ are explicit constants. 
Latter, T\'oth and Zhai \cite{Zhai_2} improved the error term in \eqref{WGN1} to
$O(x^{3/2}(\log x)^{6.5})$
by means 
of multiple complex function theory and exponential sums theory.  

For $r=3$, by using a multidimensional Perron's formula and the complex integration method,
T\'oth  and Zhai \cite{group_3} give the following asymptotic formula
\begin{equation*}
	C_3^{(1)}(x)=x^3\left(\sum_{j=0}^{7}c_j\log^jx\right)+O(x^{\frac{8}{3}+\varepsilon}),
\end{equation*}
where $c_j$ $(0\leqslant j\leqslant7)$ are explicit constants.

Let $x\geqslant 1$ be a real number and $r\in \mathbb{N}$, define 
\begin{equation}\label{D_rs}
	D_{rs}(x):=\sum_{n_1n_2\dots n_r\leqslant x}s(n_1,\dots,n_r),
\end{equation}
\begin{equation}\label{D_rc}
	D_{rc}(x):=\sum_{n_1n_2\dots n_r\leqslant x}c(n_1,\dots,n_r).
\end{equation}
Then $D_{rs}(x)$ and $D_{rc}(x)$ represent the total number of subgroups and the total number of cyclic subgroups of  finite abelian groups with order not exceed $x$ and rank  not exceed $r$, respecitively.

For $r=2$,
Sui and Liu \cite{group_2} proved that 
\begin{equation*}
	D_{2s}(x)=xP_4(\log x)+O(x^{\frac{2}{3}}\log^6x),\qquad D_{2c}(x)=xR_4(\log x)+O(x^{\frac{2}{3}}\log^6x),
\end{equation*}
where $P_4(u)$ and $R_4(u)$ are polynomials in $u$ of degree 4 with computable coefficients.  
Let 
$$\Delta_{2s}(x):=D_{2s}(x)-xP_4(\log x), \qquad\qquad 
\Delta_{2c}(x):=D_{2c}(x)-xR_4(\log x). 
$$
They also give the upper bounds of the integral mean-square of $\Delta_{2s}(x)$ and $\Delta_{2c}(x)$ as 
\begin{equation*}
\int_{1}^{T}\Delta_{2s}^2(x){\rm{d}}s\ll T^{\frac{77}{36}+\varepsilon},
\qquad\qquad
\int_{1}^{T}\Delta_{2c}^2(x){\rm{d}}c\ll T^{\frac{77}{36}+\varepsilon}.
\end{equation*}
Kiuchi and  Saad Eddin \cite{weighted}		discussed the weighted average number of subgroups of $\mathbb {Z} _m\times\mathbb {Z} _n $ with $mn\leqslant x$. In particular, for 
\begin{equation*} 
	\tilde{D}_{2s}(x):=\sum_{mn \leqslant x}s(m, n)\log\frac{x}{mn},
\qquad\qquad
	\tilde{D}_{2c}(x):=\sum_{mn \leqslant x}c(m,n)\log\frac{x}{mn},
\end{equation*}
they give the following asymptotic formulas
\begin{equation*}
	\tilde{D}_{2s}(x)=x\tilde{P}_4(\log x)+O(x^{\frac{2}{3}}\log^6x),
	\qquad\quad \tilde{D}_{2c}(x)=x\tilde{R}_4(\log x)+O(x^{\frac{2}{3}}\log^6x),
\end{equation*}
where $\tilde{P}_4(u)$ and $\tilde{R}_4(u)$ are polynomials in $u$ of degree 4 with computable coefficients.

In the case of  $r=3$ in \eqref{D_rc},  $D_{3c}(x)$  represent the total number of cyclic subgroups of  finite abelian groups with order not exceed $x$ and rank  not exceed $3$.
In this paper we shall  give an asymptotic formula for $D_{3c}(x)$, 
and get an upper bound for the integral mean-square of the error term in that asymptotic formula. 
In particular, We will prove  the following theorems, by using  Perron's formula and the 
complex integration method.


\begin{theorem1}\label{Theorem1}
	With the notation in \eqref{D_rc},
	$D_{3c}(x)$  represent the total number of cyclic subgroups of  finite abelian groups with order not exceed $x$ and rank  not exceed $3$, 
	we have
	\begin{equation*}
		D_{3c}(x)=xP_9(\log x)+O_{\varepsilon}(x^{\frac{5}{6}+\varepsilon }),
	\end{equation*}
	where $P_9(u)$ is a polynomial in $u$ of degree $9$ with computable coefficients. 
\end{theorem1}

\begin{theorem2}\label{Theorem2}
	With the notation from Theorem 1.1, and
	\[\Delta_{3c}(x):=D_{3c}(x)-xP_9(\log x),\]  
	we have
	\begin{equation*}
		\int_{1}^{T}\Delta_{3c}^2(x){\rm{d}}x\ll_{\varepsilon}T^{\frac{13}{5}+\varepsilon}
	\end{equation*}
	for $T\geqslant 10$. 
\end{theorem2}
Note that, our proof of Theorem \ref{Theorem2} depend on T\'oth and Zhai \cite[Proposition 3.1]{group_3}, which requires 
$\sigma>\frac{4}{5}$, and we have taken $\sigma=\frac{4}{5}+\varepsilon$ in our proof. So that $T^{\frac{13}{5}+\varepsilon}$ is the
best upper bound we can get in Theorem \ref{Theorem2} using \cite[Proposition 3.1]{group_3}.


Our proof of  Theorem \ref{Theorem1} depends on the Proposition 3.1 in T\'oth and Zhai \cite{group_3}, whose proof needs the following formula
from \cite[Theorem 1]{TL}
\begin{equation*}
	c(n_1,\dots,n_r)=\sum_{d_1|n_1,\dots,d_r|n_r}\dfrac{\varphi(d_1)\dots\varphi(d_r)}{\varphi([d_1,\dots,d_r])}.
\end{equation*}
However, there is no such a closed formula for $s(n_1,\dots,n_r)$ for $r\geqslant 3$ yet.  
Hence we still need to find a new way to discuss $D_{rs}(x)$. 

In the following section, $s$ will always denote a complex number with real part $\mathscr{R}s$.

\section{Preliminary lemmas}
In this section, we will give some lemmas, which will be needed in our proofs.

We need the following Proposition to write the triple Dirichlet series with coefficients $c(n_1,n_2,n_3)$ to a product of zeta-functions and a triple Dirichlet series. 

\begin{lemma2}\label{Lemma2}\cite[Proposition 3.1]{group_3}
	Let 
	\begin{equation*}
		C(s_1,s_2,s_3):=\sum_{n_1=1}^{+\infty}\sum_{n_2=1}^{+\infty}\sum_{n_3=1}^{+\infty}
		\frac{c(n_1,n_2,n_3)}{n_1^{s_1}n_2^{s_2}n_3^{s_3}}.
	\end{equation*}
	For $j=1,2,3$,  write   $s_j=\sigma_j+it_j$ with $\sigma_j, t_j\in\mathbb{R}$ and $\sigma_j>1$. Then 
	\begin{flalign*}
		C(s_1,s_2,s_3)&=\zeta^2(s_1)\zeta^2(s_2)\zeta^2(s_3)\zeta(s_1+s_2-1)\zeta(s_1+s_3-1)  
		\\
		&\ \ \  \ \times\zeta(s_2+s_3-1)\zeta(s_1+s_2+s_3-2)H(s_1,s_2,s_3),
	\end{flalign*}
	where $H(s_1,s_2,s_3)$ can be written as a triple Dirichlet series, which is absolutely convergent
	when $\sigma_1$, $\sigma_2$, $\sigma_3$ satisfy the following conditions:
	\begin{flalign*}
		&\sigma_j>\frac{1}{3}\ (j=1,2,3), \hskip33.6mm
		\sigma_1+\sigma_2+\sigma_3+\sigma_j>3\ (j=1,2,3),\ \ \ \ \  \\
		&	2\sigma_j+\sigma_i>2\ (1\leqslant j\neq i \leqslant3), \qquad\qquad
		2(\sigma_1+\sigma_2+\sigma_3)-\sigma_j>4\ (j=1,2,3). 
	\end{flalign*}
\end{lemma2}

We will use the following two version of 	Perron's formula.

\begin{lemma12}[Perron's formula]\label{Lemma12}\cite[Chapter V, Theorem 1]{A.A.K}	
	Assume that $f(s)=\sum_{n=1}^{\infty}\frac{a_n}{n^s}$  converges absolutely 	for $\sigma>1$, where $a_n\in\mathbb{C}$, $|a_n|\leqslant A(n)$, $A(n)>0$ is a monotonically increasing function,
	and 
	\begin{equation*}
		\sum_{n=1}^{\infty}|a_n|n^{-\sigma}=O((\sigma-1)^{-\alpha}),\quad \alpha>0,
	\end{equation*}
	as $\sigma\to1^{+}$. 
	Then, for any $b_0\geqslant b>1$ and $x=N+\frac12$, we have 
	\begin{equation*}
		\sum_{n\leqslant x}a_n=\dfrac{1}{2\pi i}\int_{b-iT}^{b+iT}f(s)\dfrac{x^s}{s}{\rm{d}}s
		+O\Big(\dfrac{x^b}{T(b-1)^\alpha}\Big)
		+O\Big(\dfrac{xA(2x)\log x}{T}\Big),
	\end{equation*}
	where the constants in the $O$ symbols depend only on $b_0$.
	
\end{lemma12}

\begin{lemma13}[Perron's formula]
	\label{Lemma13}
	\cite[Chapter II.2 Theorem 2.1]{TLBM}
	Let 
	\[
	F(s):=\sum_{n\geqslant1}\frac{a_n}{n^s}
	\]
	be a Dirichlet series with abscissa of convergence $\sigma_c$. 
	Put $a_x=0$ for
	$x\in\mathbb{R}\setminus\mathbb{N}^{*}$.
	Then for $\kappa>\max(0,\sigma_c)$ we have 
	\begin{equation*}
		A^{*}(x)
		:=\sum_{n<x}a_n+\frac{1}{2}a_x 
		=\frac{1}{2\pi i}\int_{\kappa-i\infty}^{\kappa+i\infty}F(s)x^s\frac{{\rm d}s}{s}\quad(x>0),
	\end{equation*}
	where the integral is conditionally convergent for
	$x\in\mathbb{R}\setminus\mathbb{N}$ and converges in the sense of Cauchy's principal value when $x\in \mathbb{N}$.
\end{lemma13}	

We need the following estimation on zeta-function.
\begin{lemma3}\label{Lemma3}\cite[\S7.5 Theorem 2]{Pan_and_Pan}
	For $l=0$ or $l=1$, $\sigma\geqslant1$ and $|t|\geqslant2$, we have 
	\begin{equation*}
		\zeta^{(l)}(\sigma+it)\ll {\rm{min}}\bigg(\frac{\sigma}{(\sigma-1)^{l+1}},(\log (|t|))^{l+1}\bigg).
	\end{equation*}
\end{lemma3}

\begin{lemma4}\label{Lemma4}\cite[Lemma 2.4]{group_3}
	Suppose $l\geqslant 0$ is a fixed integer. Then we have 
	\begin{flalign*}
		\zeta^{(l)}(\sigma+it)\ll
		\begin{cases}
			(|t|+2)^{\frac{1-\sigma}{3}}\log^{1+l}(|t|+2),\quad \frac{1}{2}\leqslant\sigma\leqslant1,\\
			(|t|+2)^{\frac{3-4\sigma}{6}}\log^{1+l}(|t|+2),\quad 0\leqslant\sigma\leqslant\frac{1}{2}.
		\end{cases}
	\end{flalign*}
\end{lemma4}

\begin{lemma10}\label{Lemma10}\cite[Theorem 8.4]{A. Ivic}
	For each fixed $\frac12<\sigma<1$, let $m(\sigma)$ be the supremum of   $m$ such that
	$\int_{1}^{T}|\zeta(\sigma+it)|^mdt\ll_{\varepsilon}T^{1+\varepsilon}$ for any $\varepsilon>0$.
	Then 
	\begin{alignat*}{2}\label{m}
		&m(\sigma)\leqslant\frac{4}{3-4\sigma},&
		&\text{\quad for\quad}\frac{1}{2}\leqslant\sigma\leqslant\frac{5}{8},\\
		&m(\sigma)\leqslant\frac{19}{6-6\sigma},&
		&\text{\quad for\quad}\frac{35}{54}\leqslant\sigma\leqslant\frac{41}{60},\\
		&m(\sigma)\leqslant\frac{12408}{4537-4890\sigma},
		&		&\text{\quad for\quad}\frac{3}{4}\leqslant\sigma\leqslant\frac{5}{6}. \\
	\end{alignat*}		
\end{lemma10}

We also need the following 	lemmas.
\begin{lemma1}\label{Lemma1}\cite[Theorem 1.3]{A. Ivic}
	Let 
	\begin{equation*}
		\zeta(s)=\frac{1}{s-1}+\sum_{k=0}^{+\infty}\frac{(-1)^k\gamma_k}{k!}{(s-1)}^k
	\end{equation*}
	be the Laurent expansion of $\zeta(s)$ in the neighborhood of its pole $s=1$. 
	Then 
	\begin{equation*}
		\gamma_k=\lim\limits_{x\to\infty}
		\bigg(\sum_{n\leqslant x}\dfrac{(\log n)^k}{n}-\dfrac{(\log x)^{k+1}}{k+1}\bigg),
	\end{equation*}
	and, 
	in particular,
	\begin{equation*}
		\gamma_0
		=\lim\limits_{N\to\infty}\bigg(1+\frac12+\dots+\frac1N-\log N\bigg)
		=0.5772157\dots
	\end{equation*}
	is Euler's constant.
\end{lemma1}

\begin{lemma5}\label{Lemma5}\cite[Chapter 8 (8.73)]{A. Ivic}
	For each real $\sigma$, let $\eta(\sigma)$ be the infimum of $a\geqslant 0$
	such that 
	$\zeta(\sigma+it)\ll t^{a}$. Then 
	\begin{equation*}
		\eta(\sigma)\leqslant \frac{1-\sigma}{5} \quad \text{for}\quad\frac56\leqslant\sigma\leqslant1.
	\end{equation*}
\end{lemma5}

\begin{lemma9}\label{Lemma9}\cite[Theorem 8.3]{A. Ivic}
	For any fixed $A\geqslant4$, define $M(A)$ to be the infimum of all $M$ such that 
	$\int_{1}^{T}|\zeta(\frac{1}{2}+it)|^Adt\ll_{\varepsilon}T^{M+\varepsilon}$ for
	any $\varepsilon>0$. Then	
	\begin{equation*} 
		M(A)\leqslant 
		\begin{cases} 
			\setstretch{1.5} 1+(A-4)/8,  &4\leqslant A\leqslant12,\\ 
			\setstretch{1.5} 2+3(A-12)/22, &12\leqslant A\leqslant\frac{178}{13},\\ \setstretch{1.5} 1+35(A-6)/216, &A\geqslant\frac{178}{13}. 
		\end{cases} 
	\end{equation*}
\end{lemma9}

\begin{lemma7}\label{Lemma7}\cite[Theorem 1.6]{A. Ivic}
	For any $s\in \mathbb{C}$, we have 
	\begin{equation*}
		\zeta(s)=\chi(s)\zeta(1-s),
	\end{equation*}
	where 
	\begin{equation*}
		\chi(s):=\dfrac{(2\pi)^s}{2\Gamma(s)\cos \frac{\pi s}{2}}
		\ll
		(2\pi)^{\sigma-\frac{1}{2}}|t|^{\frac{1}{2}-\sigma}
		\ \ (0<\sigma<1).
	\end{equation*}
\end{lemma7}


\begin{lemma11}[Plancherel's Theorem]\label{Lemma11}
	Assume $u\in L^1(\mathbb{R}^n)\cap L^2(\mathbb{R}^n)$. Then $\hat{u}$, $\check{u}\in L^2(\mathbb{R}^n)$ and
	\begin{equation*}
		\Vert\hat{u}\Vert_{L^2(\mathbb{R}^n)}
		=\Vert\check{u}\Vert_{L^2(\mathbb{R}^n)}
		=\Vert u\Vert_{L^2(\mathbb{R}^n)},
	\end{equation*}
	where  $\hat{u}$ is the Fourier transform of $u$ defined by 
	\begin{equation*}
		\hat{u}(y):=\dfrac{1}{(2\pi)^{n/2}}\int_{\mathbb{R}^n}e^{-ixy}u(x){\rm{d}}x\quad (y\in\mathbb{R}^n)
	\end{equation*} 
	and $\check{u}$ is the inverse Fourier transform of $u$ defined  by 
	\begin{equation*}
		\check{u}(y):=\dfrac{1}{(2\pi)^{n/2}}\int_{\mathbb{R}^n}e^{ixy}u(x){\rm{d}}x\quad (y\in{\mathbb{R}}^n).
	\end{equation*}
	Since $\vert e^{\pm ixy}\vert=1$ and $u\in L^1(\mathbb{R}^n)$, these integrals converge for each
	$y\in\mathbb{R}^n$.
\end{lemma11}

\begin{lemma8}[H\"older Inequation]\label{Holder}
	Suposse $f_1(x), \dots, f_n(x)\in C[a,b]$ and $ \lambda_1, \dots, \lambda_n>0$ 
	with $\frac{1}{\lambda_1}+\cdots+\frac{1}{\lambda_n}=1$. Then
	\begin{equation*}
		\int_{a}^{b}|f_1(x)\cdots f_n(x)|{\rm{d}}x\leqslant
		\bigg(\int_{a}^{b}|f_1(x)|^{\lambda_1}{\rm{d}}x \bigg)^{\frac{1}{\lambda_1}}\cdots
		\left(\int_{a}^{b}|f_n(x)|^{\lambda_n}{\rm{d}}x \right)^{\frac{1}{\lambda_n}}.
	\end{equation*}
\end{lemma8}

\begin{lemma6}\label{Lemma6}\cite[Lemma 8.3]{A. Ivic}
	Let $F(s)$ be regular in the region $\mathscr{D}:\alpha\leqslant\sigma\leqslant\beta$, $t\geqslant1$ 
	and $F(s)\ll e^{Ct^2}$ for 
	$s\in \mathscr{D}$. Then, for any fixed $q>0$ and $\alpha<\gamma<\beta$, we have
	\begin{equation*}
		\int_{2}^{T}|F(\gamma+it)|^q{\rm{d}}t\ll
		{\left(\int_{1}^{2T}|F(\alpha+it)|^q{\rm{d}}t+1 \right)}^{\frac{\beta-\gamma}{\beta-\alpha}}
		{\left(\int_{1}^{2T}|F(\beta+it)|^q{\rm{d}}t+1 \right)}^{\frac{\gamma-\alpha}{\beta-\alpha}}.
	\end{equation*}
\end{lemma6}

\section{Proof of Therorem 1.1}

\subsection{Applying Perron's Formula}

By Lemma \ref{Lemma2}, we have 
\begin{equation}\label{el1}
\sum_{l=1}^{+\infty}\sum_{m=1}^{+\infty}\sum_{n=1}^{+\infty}	\frac{c(l,m,n)}{(lmn)^{s}}
	=
	\zeta^6(s)\zeta^3(2s-1)\zeta(3s-2)T(s),
\end{equation}
where $T(s)$ can be expressed as a triple Dirichlet series, and it is uniformly 
bounded with respect to the imaginary part for $\mathscr{R}s>\frac{4}{5}$, i.e.
\begin{equation}\label{T(s)}
	T(s)\ll 1\quad(\mathscr{R}s>\frac{4}{5}).
\end{equation}
Write $c_i :=T^{(i)}(1)$ and  $T(s)=\underset{n\geqslant1}{\sum}\frac{t(n)}{n^s}$.
Since $T(s)$ is analytic for $s=\frac{4}{5}+\frac{\varepsilon}{2}$, 
then we have 
\begin{equation}\label{el2}
t(n)\ll_{\varepsilon}n^{\frac{4}{5}+\frac{\varepsilon}{2}}.
\end{equation}

Define
\begin{equation*}
f(k):=\underset{k=lmn}{\sum}c(l,m,n).
\end{equation*}
Then 
$D_{3c}(x)=\underset{k\leqslant x}{\sum}f(k)$. 
By formula \eqref{el1} and \eqref{el2}, we have 
\begin{flalign}\label{el3}
f(k)
&=\sum_{n_1n_2^2n_3^2n_4^2n_5^3n_6=k}\tau_6(n_1)n_2n_3n_4n_5^2t(n_6) \notag 
\\
&\ll_{\varepsilon}\sum_{n_1n_2^2n_3^2n_4^2n_5^3n_6=k}n_1^{\varepsilon}n_2n_3n_4n_5^2n_6^{\frac{4}{5}+\frac{\varepsilon}{2}} \notag 
\\
&\ll_{\varepsilon}k^{\frac{4}{5}+\frac{2\varepsilon}{3}}\sum_{n_1n_2^2n_3^2n_4^2n_5^3n_6=k} 1
\\
&\ll_{\varepsilon}k^{\frac{4}{5}+\varepsilon}.  \notag
\end{flalign}
Applying  the Perron's formula in  Lemma \ref{Lemma12} with $A(n)=n^{\frac{4}{5}+\varepsilon}$, 
$b=1+\frac{1}{\log x}$, $\alpha=10$ and 
 $\sqrt{x}\leqslant T\leqslant x^2$, 
we get 
\begin{equation}\label{el4}
D_{3c}(x)
=I(x,T)
  +O_{\varepsilon}\bigg(\frac{x^{\frac{9}{5}+\varepsilon}}{T}\bigg),
\end{equation}
where 
\begin{equation*}
I(x,T)
:=\frac{1}{2\pi i}\int_{b-iT}^{b+iT}\zeta^6(s)\zeta^3(2s-1)\zeta(3s-2)T(s)x^ss^{-1}{\rm{d}}s.
\end{equation*}

\subsection{Evaluation of $I(x,T)$}
In this subsection we will evaluate the integral $I(x,T)$. 

Consider the rectangle domain formed by the four points
$s=b\pm iT$ and $s=\frac{5}{6}\pm iT$. 
From \eqref{el1} we see that, in this domain,  
\begin{equation}\label{g(s)}
g(s):=\zeta^6(s)\zeta^3(2s-1)\zeta(3s-2)
\end{equation}
has only one simple pole, which is $s=1$. 
From the residue theorem we get 
\begin{equation}\label{el5}
I(x,T)=J(x,T)+H_1(x,T)+H_2(x,T)-H_3(x.T),
\end{equation}
where 
\begin{align*}
&J(x,T):=\underset{s=1}{\rm{Res}}\ g(s)T(s)\frac{x^s}{s},
\\
&H_1(x,T):=\frac{1}{2\pi i}\int_{\frac{5}{6}+iT}^{b+iT}g(s)T(s)\frac{x^s}{s}{\rm{d}}s,
\\
&H_2(x,T):=\frac{1}{2\pi i}\int_{\frac{5}{6}-iT}^{\frac{5}{6}+iT}g(s)T(s)\frac{x^s}{s}{\rm{d}}s,
\\
&H_3(x,T):=\frac{1}{2\pi i}\int_{\frac{5}{6}-iT}^{b-iT}g(s)T(s)\frac{x^s}{s}{\rm{d}}s.
\end{align*}

\subsubsection{Evaluation of $J(x,T)$}
Using  Laurent theorem, 
the complex function $T(s)$ can be uniquely written as 
\begin{equation*}
T(s)=\sum_{n=0}^{+\infty}c_n(s-1)^n,
\end{equation*}
where 
\begin{equation*}
	c_n
	:=
	  \frac{1}{2\pi i}
	   \int_{\Gamma}\dfrac{T(s)}{(s-1)^{n+1}}{\rm{d}}s
\end{equation*}
and $\Gamma$ is a circle centred at $s=1$. 
Write
\begin{equation*}
T^*(s)
: =
\sum_{k=0}^{N}\frac{T^{(k)}(1)}{k!}(s-1)^k.
\end{equation*}
From Lemma \ref{Lemma1} we see that
\begin{equation*}
\underset{s=1}{\rm{Res}}\ g(s)\frac{T(s)x^s}{s}=\underset{s=1}{\rm{Res}}\ g(s)\frac{T^*(s)x^s}{s}.
\end{equation*}
By numerically calculating we have 
\begin{equation}\label{el6}
J(x,T)=\underset{s=1}{\rm{Res}}\ \frac{T^*(s)x^s}{s}
=xP_9(\log x), 
\end{equation}
where $P_9(u)=\sum_{r=0}^{9}A_r u^r$ is a polynomial in $u$ with 
$$A_9=\dfrac{T(1)}{8709120},
\hskip 20mm 
A_8=\dfrac{T'(1)+T(1)(15\gamma  -1)}{967680}.
$$ 
The expression of $A_7\sim A_0$ are complicate so that we omit them here.

\subsubsection{Evaluation of $H_j(x,T)$ for $j=1,2,3$} 
We estimate $H_1(x,T)$ first. 

Let $G(s):=g(s)T(s)\frac{x^s}{s}$ and $s=\sigma+it$.
For $1\leqslant\sigma<1+\frac{1}{\log x}$, by using Lemma \ref{Lemma3}, we have
$$
|\zeta(s)| \ll\log T
$$
for $s=\sigma+iT, 2\sigma-1+2iT, 3\sigma-2+3iT$ and $T>10$.
For $\frac56 \leqslant \sigma \leqslant 1$, 
notice $2\sigma-1\in(\frac{2}{3},1)$,  $3\sigma-2\in(\frac{1}{2},1)$, using Lemma \ref{Lemma4} and Lemma 
\ref{Lemma5}, we get
$	|\zeta(\sigma+iT)|^6\ll T^{\frac{6(1-\sigma)}{5}}$,
$	|\zeta(2\sigma-1+2iT)|^3\ll T^{2(1-\sigma)}\log^3T$,
$	|\zeta(3\sigma-2+3iT)|\ll T^{1-\sigma}\log T$. 
Then, using \eqref{T(s)}, we have 
\begin{flalign*}
	G(s)   &\ll|\zeta(\sigma+iT)|^6|\zeta(2\sigma-1+2iT)|^3|\zeta(3\sigma-2+3iT)|T^{-1} \notag \\ 
	&\ll
	\begin{cases}
		x^{\sigma}T^{\frac{16}{5}-\frac{21}{5}\sigma}\log^4T,\ &\frac{5}{6}\leqslant\sigma\leqslant1; \\
		xT^{-1}\log ^{10}T, \  &\sigma>1.
	\end{cases}
\end{flalign*}
Notice that $x^\frac12\leqslant T\leqslant x^2$, we have 
\begin{flalign}\label{el8}
H_1(x,T)&\ll_{\varepsilon}\int_{\frac{5}{6}}^{1}x^\sigma 
T^{\frac{16}{5}-\frac{21}{5}\sigma}\log^4 T{\rm{d}}\sigma+
\int_{1}^{1+\frac{1}{\log x}}xT^{-1}\log^{10}T{\rm d}\sigma \notag \\
&\ll_{\varepsilon}T^{\frac{16}{5}}\log^4T\left(\dfrac{x}{T^{\frac{21}{5}}} \right)^{\sigma}\Big/
\log \left(\frac{x}{T^{\frac{21}{5}}}\right)\Big|_{\frac{5}{6}}^{1}+\dfrac{x}{\log x}T^{-1}\log^{10}T   \\
&\ll_{\varepsilon}x^{\frac{41}{60}}\log^3T.  \notag
\end{flalign}

Similarly, we can get
\begin{equation}\label{el9}
H_3(x,T)\ll_{\varepsilon}x^{\frac{41}{60}}\log^3T.
\end{equation}

Now we estimate $H_2(x,T)$. It is not difficult to see that 
\begin{flalign*}
H_2(x,T)
&\ll_{\varepsilon}\int_{-T}^{T}\dfrac{x^{\frac{5}{6}}|\zeta({\frac{5}{6}+it})|^6|\zeta(\frac{2}{3}+2it)|^3
|\zeta(\frac{1}{2}+3it)|}{|t|+1}{\rm{d}}t \\
&\ll_{\varepsilon}x^{\frac{5}{6}}
 \bigg(1+\int_{1}^{T}\dfrac{|\zeta({\frac{5}{6}+it})|^6|\zeta(\frac{2}{3}+2it)|^3
|\zeta(\frac{1}{2}+3it)|}{t}{\rm{d}}t\bigg) \\
&\ll_{\varepsilon}x^{\frac{5}{6}}+x^{\frac{5}{6}}\int_{1}^{T}\dfrac{S_1^6S_2^3S_3}{t}{\rm{d}}t,
\end{flalign*}
where 
\begin{equation*}
	S_1=\left|\zeta\bigg(\frac{5}{6}+it\bigg)\right|,\qquad
	S_2=\left|\zeta\bigg(\frac{2}{3}+2it\bigg)\right|,\qquad
	S_3=\left|\zeta\bigg(\frac{1}{2}+3it)\right|.
\end{equation*}
Applying the H\"older Inequation in Lemma \ref{Holder}, we  get
\begin{equation*}
\int_{1}^{T}\dfrac{S_1^6S_2^3S_3}{t}{\rm{d}}t\ll_{\varepsilon}
\left(\int_{1}^{T}\dfrac{S_1^{24}}{t}{\rm{d}}t \right)^{\frac{1}{4}}
\left(\int_{1}^{T}\dfrac{S_2^{12}}{t}{\rm{d}}t \right)^{\frac{1}{4}}
\left(\int_{1}^{T}\dfrac{S_3^{2}}{t}{\rm{d}}t \right)^{\frac{1}{2}}.
\end{equation*}
Here we have 
$24\leqslant\frac{12408}{4537-4890\times\frac{5}{6}}$, 
$12\leqslant\frac{19}{6-6\times\frac{2}{3}}$ 
and 
$2\leqslant\frac{4}{3-4\times\frac{1}{2}}$, so that we can apply 
Lemma \ref{Lemma10}
to get
\begin{equation}\label{el10}
H_2(x,T)\ll_{\varepsilon}x^{\frac{5}{6}+\varepsilon}.
\end{equation}

\subsubsection{Completion of the proof of Theorem \ref{Theorem1}} 
Now we can evaluate $D_{3c}(x)$. 

Take $T=x$. 
Then,   \eqref{el4}--\eqref{el10} imply
\begin{equation*}
	D_{3c}(x)=xP_9(\log x)+O_{\varepsilon}(x^{\frac{5}{6}+\varepsilon}),
\end{equation*}             
 where $P_9(u)$ is the polynomial defined in \eqref{el6}.
 Then we complete the proof of Theorem \ref{Theorem1}.

\section{Proof of Therorem 1.2}

In this section, we shall employ Fourier transform and the Plancherel's Theorem to give an upper bound of the
mean-square estimate of $\Delta_{3c}(x)$.

\subsection{Fourier transform}
\indent
Consider the rectangle domain formed by $s=\sigma $ and $s=2$. By the Perron's formula in Lemma \ref{Lemma13} we get
\begin{flalign*}
	\prescript{ }{ }{\sum _{n< x}}^{'}f(n)&=
	\dfrac{1}{2\pi i}\int_{2-i\infty}^{2+i\infty}\zeta^6(s)\zeta^3(2s-1)\zeta(3s-2)T(s)x^ss^{-1}{\rm{d}}s,
\end{flalign*}
where $\prescript{ }{ }{\sum^{'} _{n< x}}f(n)
:=\prescript{ }{ }{\sum _{n< x}} f(n)+\frac{1}{2}f(x)$ for  $x\in\mathbb{Z}$ 
and 
$\prescript{ }{ }{\sum^{'} _{n< x}}f(n)
:=\prescript{ }{ }{\sum _{n< x}} f(n)$ 
for $x\notin\mathbb{Z}$.
Then applying the residue theorem and \eqref{el6}, we have 
\begin{flalign}\label{el11}
	\sum_{k< x}f(k)
	&=xP_9(\log x)+\frac{1}{2\pi}\int_{-\infty}^{+\infty}
	h(t)x^{\sigma +it}{\rm{d}}t+O_{\varepsilon}(x^{\frac{4}{5}+\varepsilon}),
\end{flalign}
where 
\begin{equation*}
	h(t):=
	\zeta^6(\sigma +it)\zeta^3(2\sigma -1+2it)\zeta(3\sigma -2+3it)T(\sigma +it)(\sigma +it)^{-1}.
\end{equation*}
Write
\begin{equation*}
	\delta_{3c }(x)=\frac{1}{2\pi}\int_{-\infty}^{+\infty}h(t)x^{\sigma +it}{\rm d}t.
\end{equation*}
Then $\Delta_{3c}(x) = \delta_{3c }(x) + O_{\varepsilon}(x^{\frac{4}{5}+\varepsilon})$, and $\delta_{3c }(x)$ is the Fourier inversion of $h(t)$. So, we can write $h(t)$ as 
\begin{equation*}
	h(t)=\frac{1}{(2\pi)^{1/2}}\int_{-\infty}^{+\infty}(2\pi)^{1/2}
	\delta_{3c}(x)x^{-\sigma }e^{-it\log x}{\rm{d}}\log x.
\end{equation*}
Applying the Plancherel's Theorem in Lemma \ref{Lemma11}, we have
\begin{flalign}\label{Plancherel}
	\int_{-\infty}^{+\infty}|h(t)|^2{\rm d}t
	=2\pi\int_{0}^{+\infty}|\delta_{3c}(x)|^2x^{-2\sigma -1}{\rm d}x.
\end{flalign}

\subsection{Completion of the proof of Theorem \ref{Theorem2}}
Take $\frac{4}{5}<\sigma <1$. 
For any $U\geqslant2$, using Lemma \ref{Lemma7} and \eqref{T(s)}, we have 
\begin{flalign}\label{el13}
	&\int_{U}^{2U} |h(t)|^2 {\rm{d}}t\notag \\
&
	\ll\int_{U}^{2U}|\zeta(\sigma +it)|^{12}|\zeta(2\sigma -1+2it)|^6|\zeta(3-3\sigma -3it)|^2U^{2(\frac{5}{2}-3\sigma )}U^{-2}{\rm{d}}t \\
	&
	\ll U^{3-6\sigma }K(U),\notag 
\end{flalign}
where 
\begin{equation*}
	K(U):=\int_{U}^{2U}W_1^{12}W_2^6W_3^2{\rm{d}}t,
\end{equation*}
\begin{equation*}
	W_1:=|\zeta(\sigma +it)|,\qquad W_2:=|\zeta(2\sigma -1+2it)|,\qquad W_3:=|\zeta(3-3\sigma -3it)|.
\end{equation*}
Applying the H\"older Inequation in Lemma \ref{Holder} with $\lambda_1=\frac{33}{20}$, $\lambda_2=\frac{66}{15}$ and $\lambda_3=6$,
we  get
\begin{equation}\label{el14}
	K(U)
	\ll
	\left(\int_{U}^{2U}W_1^{\frac{99}{5}}{\rm{d}}t \right)^{\frac{20}{33}}
	\left(\int_{U}^{2U}W_2^{\frac{396}{13}}{\rm{d}}t \right)^{\frac{13}{66}}
	\left(\int_{U}^{2U}W_3^{12}{\rm{d}}t \right)^{\frac{1}{6}}.
\end{equation}
Then,  using Lemma \ref{Lemma10}, Lemma \ref{Lemma6} and Lemma \ref{Lemma9}, we have
$$
	\int_{U}^{2U}W_1^{\frac{99}{5}}{\rm{d}}t
	\ll
	U^{1+\varepsilon},
$$
\begin{flalign*}
	&  \int_{U}^{2U}W_2^{\frac{396}{13}}{\rm{d}}t
	\\
	&\ll
	\left(\int_{1}^{2U}	\left|\zeta\bigg(\frac{3}{5}+2it\bigg)\right|^{\frac{396}{13}}{\rm{d}}t +1\right)^{\frac{5/6-2\sigma  +1}{5/6-3/5}} 
	\left(\int_{1}^{2U}\left|\zeta\bigg(\frac{5}{6}+2it\bigg)\right|^{\frac{396}{13}}{\rm{d}}t +1\right)^{\frac{2\sigma -1-3/5}{5/6-3/5}} \\
	&\ll U^{1+\varepsilon}
\end{flalign*}
and
\begin{flalign*}
	&\int_{U}^{2U}W_3^{\frac{1}{3}}{\rm{d}}t
	\\
	&\ll
	\left(\int_{1}^{2U}\left|\zeta\bigg(\frac{1}{2}-3it\bigg)\right|^{\frac{1}{3}}{\rm{d}}t +1\right)^{\frac{4/19-3+3\sigma }{4/19-1/2}}
	\left(\int_{1}^{2U}\left|\zeta\bigg(\frac{4}{19}-3it\bigg)\right|^{\frac{1}{3}} {\rm{d}}t +1\right)^{\frac{3-3\sigma -1/2}{4/19-1/2}}\\
	&\ll U^{2+\varepsilon},
\end{flalign*}
So, \eqref{el14} becomes
$$
	K(U)
\ll
\left(U^{1+\varepsilon} \right)^{\frac{20}{33}}
\left( U^{1+\varepsilon} \right)^{\frac{13}{66}}
\left(U^{2+\varepsilon} \right)^{\frac{1}{6}}
\ll U^{\frac{7}{6}+\varepsilon},
$$
which implies 
\begin{flalign*}
	\int_{U}^{2U}|h(t)|^2{\rm d }t
	\ll
	U^{-\frac{19}{30}+\varepsilon}.
\end{flalign*}
Then from \eqref{Plancherel} and \eqref{T(s)} we have 
\begin{flalign}\label{el15}
	&\int_{0}^{+\infty}|\delta_{3c}(x)|^2x^{-2\sigma -1}{\rm d}x\notag 
\\
	&\ll \sum_{j=0}^{+\infty}\int_{2^j}^{2^{j+1}}|\zeta(\sigma +it)|^{12}|\zeta(2\sigma -1+2it)|^6
	|\zeta(3\sigma -2+3it)|^2|\sigma +it|^{-2}{\rm{d}}t \\
	&\ll\sum_{j=0}^{+\infty}2^{(-\frac{19}{30}+\varepsilon)j} 
	\ll1.\notag
\end{flalign}
In particular
\begin{equation*}\label{T}
	\int_{V}^{2V}\delta_{3c}^2(x){\rm{d}}x\ll V^{2\sigma +1}
\end{equation*}
for any $V \geqslant 0$. For any $T>10$, take $\eta=[\log_2T]$. Then $T\leqslant 2^{\eta+1}$. 
Taking $\sigma =\frac{4}{5}+\varepsilon$,
we can obtain
\begin{flalign*}
	\int_{1}^{T}\Delta_{3c}^2(x){\rm d}x
	&\ll\int_{1}^{T}\left( \delta_{3c}(x)+O_{\varepsilon}(x^{\frac45+\varepsilon})\right)^2{\rm d}x \\
	&\ll\sum_{j=0}^{\eta+1}\int_{2^j}^{2^{j+1}}
	\left(\delta_{3c}^2(x)+O_{\varepsilon}(x^{\frac85+\varepsilon})\right){\rm d}x \\
	&\ll T^{\frac{13}{5}+\varepsilon}.
\end{flalign*}


\begin{thebibliography}{99}
	\bibitem{Bhowmik_Wu}
	G. Bhowmik and J. Wu, 
	Zeta function of subgroups of abelian groups and average orders. 
	J. Reine Angew. Math. 530 (2001), 1-15.
	
	\bibitem{Calugareanu.G}
	G. C\u alug\u areanu, 
	The total number of subgroups of a finite abelian group. 
	Sci. Math. Japon. 60 (2004), 157-167.
	
	\bibitem{Chinta.G_Kaplan.N_Koplewitz.S}
	G. Chinta, N. Kaplan, S. Koplewitz, 
	The cotype zeta function of $\mathbb{Z}^d$, 
	Indag. Math. (N.S.) 34 (3) (2023), 643-659.
	
	\bibitem{M_Hampejs}
	M. Hampejs, N. Holighaus, L. T\'oth and C. Wiesmeyr, Representing and counting the subgroups of the
	group $\mathbb{Z}_m\times\mathbb{Z}_n$, J. Numbers 2014, art. ID 491428, 6 pp.
	
	
	\bibitem{M_N_Huxley}
	M. N. Huxley, Exponential sums and lattice points III, 
	Proc. London Math. Soc. 87 (2003), 591-609.
	
	\bibitem{A. Ivic}
	A. Ivi\'c, The Riemann Zeta-Function. Theory and Applications, Wiley, New York, 1985.



	\bibitem{A.A.K}
	A.A. Karatsuba, Basic Analytic Number Theory, Springer-Verlag, Berlin, Heidelberg, 1993.
	
\bibitem{weighted}	
I. Kiuchi 
and S.  Saad Eddin, 
On the weighted average number of subgroups of $\mathbb 
{Z} _m\times\mathbb {Z} _n $ with $mn\leqslant x$. 
Int. J. Number Theory1 8 (9) (2022),  2005-2013.

	\bibitem{W_G_Nowak}
	W.G. Nowak and L. T\'oth, On the average number of subgroups of the group 
	$\mathbb{Z}_m\times\mathbb{Z}_n$, 
	Int. J. Number Theory 10 (2) (2014), 363-374.
	
	\bibitem{Pan_and_Pan}
	CB. Pan and CD. Pan, Foundations of Analytic Number Theory. 
	Sci. Press, Beijing, 1990 (in Chinese).
	
	\bibitem{Petrillo}
	J. Petrillo, Counting subgroups in a direct product of finite cyclic groups. 
	Coll. Math. J. 42 (2011), 215–222.
	
	\bibitem{Petrgradsky.V.M}
	V. M. Petrogradsky, Multiple zeta functions and asymptotic structure of free abelian groups of finite rank. 
	J. Pure Appl. Algebra 208 (3) (2007), 1137-1158.
	
	\bibitem{R_Schmidt}
	R. Schmidt, Subgroup Lattices of Groups, de Gruyter Exp. Math. 14, de Gruyter, Berlin, 1994.
	
	\bibitem{fourier}
	Elias M. Stein and R. Shakarchi. 
	Fourier analysis: an introduction. Vol. 1. Princeton University Press, 2011.
	
	\bibitem{group_2}
	Y. Sui and D. Liu, 
	On the error term concerning the number of subgroups of the groups $Z_m\times Z_n$ with $mn\leqslant x$. 
	J. Number Theory,  216 (2020), 264-279.
	
	\bibitem{M_Suzuki}
	M. Suzuki, 
	On the lattice of subgroups of finite groups, 
	Trans. Amer. Math. Soc. 70 (1951), 345-371.
	
	
	\bibitem{Tarnauceanu2007}
	M. T{\textup{$\breve{a}$}}rn$\breve{a}$uceanu, 
	A new method of proving some classical theorems if abelian groups. 
	Southeact Asian Bull. Math.  31 (2007), 1191-1203.
	
	\bibitem{Tarnauceanu2010}
	M. T{\textup{$\breve{a}$}}rn{\textup{$\breve{a}$}}uceanu, 
	An arithmetic method of counting the subgroups of a finite abelian group. 
	Bull. Math. Soc. Sci. Math. Roumanie (N.S.) 53 (101) (2010), 373-86.
	
	\bibitem{TLBM}
	G. Tenenbaum, 
	Introduction to analytic and probabilistic number theory. Vol. 163. American Mathematical Soc., 2015.
	
	\bibitem{TL}
	L. T\'oth, On the number of cyclic subgroups of a finite Abelian group. 
	Bull. Math. Soc. Sci. Math. Roumanie (N.S.) 55 (103) (2012), no.4, 423-428.
	
	
	\bibitem{Toth.L}
	L. T\'oth, Subgroups of finite Abelian groups having rank two via Goursat’s lemma. 
	Tatra Mt. Math. Publ. 59 (2014), 93-103.
	
	
	\bibitem{Zhai_2}
	L. T\'oth and WG. Zhai, On the error term concerning the number of subgroups of the groups $\mathbb 
	{Z} _m\times\mathbb {Z} _n $ with $m, n\leqslant x$. 
	Acta Arithmetica 183 (2018): 285-299.
	
	\bibitem{group_3}
	L. T\'oth and WG. Zhai. On the average number of cyclic subgroups of the groups $\mathbb{Z}_{n_1}
	\times\mathbb{Z}_{n_2}\times\mathbb{Z}_{n_3}$ with $n_1, n_2, n_3\leqslant x$. 
	Res.  Number Theory 6 (1) (2020): Art. 12, 33 pp.
	

	
	
	
\end{thebibliography}
\end{document}